\documentclass[12pt]{article}
\usepackage{amsfonts, amsmath, amsthm}

\newtheorem{theorem}{Theorem}[section]
\newtheorem{lemma}[theorem]{Lemma}
\newtheorem{conj}[theorem]{Conjecture}
\newtheorem{cor}[theorem]{Corollary}
\newtheorem{prop}[theorem]{Proposition}

\def\PP{\mathbb{P}}
\def\QQ{\mathbb{Q}}

\def\calO{\mathcal{O}}

\def\alg{\mathrm{alg}}
\def\an{\mathrm{an}}
\def\con{\mathrm{con}}
\def\cris{\mathrm{cris}}

\def\imm{\mathrm{imm}}
\def\sep{\mathrm{sep}}
\def\perf{\mathrm{perf}}

\def\beq{\begin{equation}}
\def\eeq{\end{equation}}
\def\map#1#2#3{#1\!:\!#2 \to #3}

\def\fp{\frac{1}{p}}

\def\GL{\Gamma^L}
\def\Gperf{\Gamma^{\perf}}
\def\Gsep{\Gamma^{\sep}}
\def\Gimm{\Gamma^{\imm}}
\def\Galg{\Gamma^{\alg}}

\def\Galgcon{\Galg_{\con}}

\def\Gancon{\Gamma_{\an,\con}}
\def\Gcon{\Gamma_{\con}}
\def\Gimmcon{\Gimm_{\con}}

\def\GLcon{\GL_{\con}}
\def\Gperfcon{\Gperf_{\con}}
\def\Gsepcon{\Gsep_{\con}}

\def\Oan{\Omega_{\an}}

\def\be{\mathbf{e}}
\def\bof{\mathbf{f}}
\def\bv{\mathbf{v}}
\def\bw{\mathbf{w}}
\def\bx{\mathbf{x}}

\DeclareMathOperator{\Gal}{Gal}
\DeclareMathOperator{\Hom}{Hom}

\DeclareMathOperator{\Spec}{Spec}

\begin{document}

\title{Descent of morphisms of overconvergent $F$-crystals}
\author{Kiran S. Kedlaya}
\date{\today}

\maketitle

\begin{abstract}
Tsuzuki has conjectured that the embedding of overconvergent
$(F, \nabla)$-crystals over $k((t))$ (for $k$ a field of characteristic
$p>0$)
into the category of convergent
$(F, \nabla)$-crystals over $k((t))$ is fully faithful.
We prove Tsuzuki's conjecture restricted
to the subcategory of potentially semistable
(or quasi-unipotent) crystals, following de~Jong's proof of a slightly
weaker result. We also prove Tsuzuki's conjecture restricted to crystals
with at most two distinct slopes.
\end{abstract}

Building on our earlier paper \cite{bib:me3},
our goal is to extend the extension theorem of de~Jong
\cite{bib:dej1}
to some morphisms between overconvergent $(F, \nabla)$-crystals over
$k((t))$, for $k$ a field of characteristic $p>0$. Tsuzuki's
``Tate type'' conjecture \cite[Conjecture~2.3.2]{bib:tsu4} states that
the embedding of the category of overconvergent $(F, \nabla)$-crystals
into the category of convergent $(F, \nabla)$-crystals is fully faithful;
he has proved this for \'etale (unit-root) crystals \cite{bib:tsu2},
and it follows readily for isoclinic crystals (those with all slopes equal).

We first prove the Tate type conjecture on the subcategory of
potentially semistable $(F, \nabla)$-crystals
(in the sense of \cite{bib:me3}) is fully faithful; de~Jong's theorem is 
the case where the crystals are semistable and 
do not have logarithmic poles. It is conjectured that in
fact every overconvergent $(F, \nabla)$-crystal is potentially semistable
\cite[Conjecture~4.12]{bib:me3}, but this is only known in general
for \'etale crystals, by a result of Tsuzuki \cite{bib:tsu1}.
Using that result, we prove the Tate type conjecture
for crystals of rank 2, and more generally for crystals 
with at most two distinct slopes each.

As applications, we answer a question of Katz by showing that an $F$-crystal
over $\Omega$ is unipotent if and only if its associated Galois representation
is trivial,
and we show that all $(F, \nabla)$-crystals arising from the crystalline
cohomology of a smooth proper variety over $k((t))$ are overconvergent and
potentially semistable.

\section{Notations}
We retain the definitions and notations of our earlier paper
\cite{bib:me3} (which are in turn modeled on \cite{bib:dej1}).
For the convenience
of the reader, we summarize these in the following table. (Note:
wherever it appears in the table, $*$ represents an unspecified
decoration.)

\begin{list}%
{}{\setlength{\itemsep}{0pt}
\setlength{\parsep}{0pt}
\setlength{\labelwidth}{0.75in}
\setlength{\leftmargin}{1in}
\setlength{\labelsep}{0.25in}}
\item[$k$]
An algebraically closed field of characteristic $p>0$
\item[$\calO$]
A finite extension of the Witt ring $W(k)$
\item[$\sigma$ on $\calO$]
An automorphism of $\calO$ lifting the absolute Frobenius on $k$
\item[$\calO_0$]
The elements of $\calO$ fixed by $\sigma$
\item[$|\cdot|$]
The valuation of $\calO$ normalized so that $|p| = p^{-1}$
\item[$K$]
The field of formal power series over $k$
\item[$K^{\perf}$]
The perfect closure of $K$
\item[$K^{\sep}$]
The separable closure of $K$
\item[$K^{\alg}$]
The algebraic closure of $K$
\item[$K^{\imm}$]
The ring of series $x = \sum_{i \in I} x_i t^i$ over $k$ with $I \subseteq \QQ$ well-ordered
\item[$\Omega_t$ (= $\Omega$)]
The power series ring $\calO[[t]]$
\item[$\Gamma$]
The $p$-adic completion of $\Omega[t^{-1}]$
\item[$\sigma$ on $\Gamma$]
An endomorphism lifting $x \mapsto x^p$ compatible 
with $\sigma$ on $\calO$
\item[$\sigma_t$]
The endomorphism of $\Gamma$ with $t \mapsto t^p$
\item[$\Gamma^L$]
The $p$-adically complete
extension of $\Gamma$ with residue field $L$
\item[$\Gamma^*$]
Equal to $\Gamma^{K^*}$ for $* \in \{\perf, \sep,
\alg, \imm\}$
\item[$\Gamma^{\alg(c)}$]
The subring of $x = \sum_i x_i t^i \in \Galg$
for which for each $n \geq 0$, there exists $r_n$
such that $|x^{\sigma^{n}}-r_n|
< p^{-cn}$
\item[$\Gamma^{*}_{\con}$]
The ring of $x = \sum_{i=-\infty}^{\infty} x_i t^i \in \Gamma^*$ with
$x_i \in \calO$ and  $\liminf_{i \to \infty} v_p(x_{-i})/i > 0$
\item[$\Gancon$]
The ring of $x = \sum_{i=-\infty}^{\infty} x_i t^i$ with
$x_i \in \calO[\fp]$, $\liminf_{i \to \infty} v_p(x_{-i})/i > 0$,
and $\liminf_{i \to +\infty} v_p(x_i)/i \geq 0$ 
\end{list}

We also recall the following basic proposition \cite[Proposition~2.7]{bib:me3},
which extends \cite[Proposition~8.1]{bib:dej1}.
\begin{prop} \label{inj1}
The following multiplication maps are injective:
\begin{eqnarray*}
\Gperfcon \otimes_{\Gcon} \Gamma &\to& \Gperf \\
\Gsepcon \otimes_{\Gcon} \Gamma &\to& \Gsep \\
\Galgcon \otimes_{\Gcon} \Gamma &\to& \Galg \\
\Gimmcon \otimes_{\Gcon} \Gamma &\to& \Gimm \\
\Gperfcon \otimes_{\Gcon} \Gsep &\to& \Galg.
\end{eqnarray*}
\end{prop}

\section{Slope filtrations}
Over $\Galg$, every $F$-crystal splits as a direct sum of trivial
crystals, but this is not true for crystals over $\Gamma$ or $\Gcon$.
In this section, we bridge the gap by
constructing two filtrations of an $F$-crystal over intermediate rings.

The following lemma
is \cite[Proposition~5.5]{bib:dej1}, but we give a simplified
approach based on the proof of \cite[Proposition~2.2.2]{bib:tsu2}.
\begin{lemma}[Descending slope filtration] \label{slopefilt}
With $\sigma$ arbitrary,
let $M$ be an $F$-crystal over $R = \Gperfcon, \Galgcon$ or $\Gimmcon$.
Then after making a suitable finite extension of $\calO$,
$M$ admits a canonical filtration
\[
0 = M_0 \subseteq M_1 \subseteq \cdots \subseteq M_n = M
\]
with $M_i/M_{i-1}$ being
isoclinic of slope $\ell_i$ and
$\ell_1 > \ell_2 > \cdots >\ell_n$.
Additionally, if $R=\Galgcon$ or $\Gimmcon$,
$M_i/M_{i-1}$ is spanned by eigenvectors of slope $\ell_i$.
\end{lemma}
\begin{proof}
The desired filtration clearly exists and is unique over $\Gperf,
\Galg, \Gimm$, and in fact all of the eigenvectors
of $M$ are defined over $\Galg$ or $\Gimm$.
It suffices to show that for $R = \Galgcon$ or $R =
\Gimmcon$, the eigenvectors of slope $\ell_1$ are defined over
$R$. (In the case $R = \Gperfcon$, the space they span is then
defined over $\Gperf \cap \Galgcon = \Gperfcon$.)

Choose $t \in \Gperfcon$ such that $t^\sigma = t$;
we will suppress $t$ in the notation of the semi-valuation $v_{t,n}$.
Choose an isobasis
$\be_1, \dots, \be_n$ of $M$ over $\Galg$ or $\Gimm$
with $F\be_i = \lambda_i \be_i$ for some $\lambda_i \in \calO_0$,
such that $|\lambda_1| \geq \cdots \geq |\lambda_n|$.
Let $k$ be the largest integer such that $|\lambda_k| =
p^{-\ell_1}$; we must show that $\be_1, \dots, \be_k$
are overconvergent.

Pick $i \in \{1, \dots, k\}$, so that $\be_i$ is an eigenvector of
maximum slope.
Now choose an isobasis
$\bof_1, \cdots, \bof_n$ of $M$ over $R$
on which $F$ acts via a matrix $A$ for which $B = \lambda_i A^{-1}$
is integral, and write $\be_i =
\sum_j c_j \bof_j$. By rescaling $\be_i$, we may assume that
the $c_j$ are integral.

Because $M$ is defined over an overconvergent ring,
we may choose $c,d >0$ such that
$v_{n}(B_{jk}) \geq -cn+d$ for all $n$.
We will show by induction on $n$ that $v_{n}(c_j) \geq -cn/(p-1)+d/(p-1)$.
Choose $j$ to minimize $v_{n}(c_j)$; the equality $F\be_i = \lambda_i \be_i$
implies the equation $c_j^\sigma = \sum_k B_{jk} c_k$ and hence the inequality
\[
p v_{n}(c_j) \geq \min_{k,m} \{ v_{m}(B_{jk}) + v_{n-m}(c_k) \}.
\]
If the minimum on the right side is achieved with $m = 0$,
we may replace the term $v_n(c_k)$ by $v_n(c_j)$ while preserving the
inequality, leading to the conclusion $(p-1)v_n(c_j) \geq -cn +d$. Otherwise,
we have by
induction
\begin{align*}
v_n(c_j) &\geq \frac 1p \left( -cm+d + \frac{-c(n-m) + d}{p-1} \right) \\
&= - \frac{cn}{p(p-1)} - \frac{cm(p-2)}{p(p-1)} + \frac{d}{p-1} \\
&\geq -\frac{cn}{p(p-1)} - \frac{cn(p-2)}{p(p-1)} + \frac{d}{p-1} \\
&= -\frac{cn}{p} + \frac{d}{p-1} \geq -\frac{cn}{p-1} + \frac{d}{p-1}.
\end{align*}
In either case, the induction follows, and we conclude that the
$c_i$ are convergent.
Therefore the submodule spanned by $\be_1, \dots, \be_i$ is indeed defined over
$R$ and so is the desired $M_1$.
\end{proof}

\begin{lemma}[Ascending slope filtration] \label{slopefilt1a}
With $\sigma$ arbitrary,
let $M$ be an $F$-crystal over $R = \Gamma, \Gsep, \Gperf, \Galg$ or $\Gimm$.
Then after a suitable finite extension of $\calO$,
$M$ admits a canonical filtration
\[
0 = M_0 \subseteq M_1 \subseteq \cdots \subseteq M_n = M
\]
with $M_i/M_{i-1}$ isoclinic of slope $\ell_i$ and
$\ell_1 < \ell_2 < \cdots <\ell_n$.
Additionally, if  $R = \Gsep, \Galg$ or $\Gimm$,
$M_i/M_{i-1}$ is spanned by eigenvectors of slope $\ell_i$.
\end{lemma}
\begin{proof}
Let $\ell_1 < \cdots < \ell_n$ be the distinct slopes of
$M$. It suffices to show that for $R = \Gsep$, the eigenvectors of $M$
over $\Galg$ of slope $\ell_1$ are defined over $\Gsep$.

Pick $\lambda \in \calO_0$ with $|\lambda| = p^{-\ell_1}$,
choose an isobasis
$\be_1, \dots, \be_n$ of $M$ over $\Gsep$ on which $F$ acts by a matrix $A$ such that
$\lambda^{-1}A$ is integral and congruent to the projection onto
the span of $\be_1, \dots, \be_k$ modulo $\pi$.
Let $\bv$ be a vector over $\Galg$ such that
$F\bv = \lambda \bv$. We show by
induction on $j$ that $\bv$ is congruent to a vector over $\Gsep$
modulo $\pi^j$, starting with $j=0$. Suppose that $\bv = \bw +
\pi^j \bx$ with $\bw$ defined over $\Gsep$. From $A\bv^\sigma =
\lambda \bv$ we get
\[
(\lambda^{-1}A \bx^\sigma - \bx) = \pi^{-j}(\bw - \lambda^{-1} A
\bw^\sigma).
\]
Reducing modulo $\pi$, we obtain equations of the form $x_i^p -
x_i = c_i$ for $1 \leq i \leq k$ and $-x_i = c_i$ for $i>k$.
Given that $c_i \in \Gsep$, we conclude that $x_i \in \Gsep$,
so that $\bv$ is congruent to a vector over $\Gsep$ modulo
$\pi^{j+1}$.

The induction complete, we conclude that $\bv$ is defined over
$\Gsep$, as desired.
\end{proof}

In passing, we note that the slope filtrations imply that
$(F, \nabla)$-crystals over $\Gcon$ may admit morphisms either over
$\Galgcon$ or $\Gsep$ that do not descend to $\Gcon$. Namely,
Lemma~\ref{slopefilt} implies that any $F$-crystal over $\Gcon$ admits
a nonzero map over $\Galgcon$ (resp.\ $\Gsep$) from a trivial crystal, 
whose image is the submodule generated by eigenvectors of maximum
slope (resp.\ minimum slope).
Thus it suffices to observe that there actually exist
$(F, \nabla)$-crystals over $\Gcon$ that have more than one distinct slope but which
are irreducible. In fact
such crystals exist over $\Omega$; the crystalline Dieudonn\'e module of a family
of elliptic curves with ordinary generic fibre and supersingular special fibre
is such a crystal. (For an example given by explicit equations, see the discussion
of the Bessel isocrystal in \cite{bib:tsu3}.)

\begin{lemma} \label{lem:filterperf}
Let $M$ be an $F$-crystal over $\Gsep$ which is not isoclinic,
and whose largest and second largest slope (not counting multiplicity)
differ by $c$. Then
the eigenvectors of $M$ of maximum slope
are defined over $\Gamma^{\alg(c)}[\fp]$.
\end{lemma}
\begin{proof}
Lemma~\ref{slopefilt1a} allows us to construct an isobasis $\be_1, \dots, \be_n$
of $M$ such that
\[
F\be_i = \lambda_i \be_i + \sum_{j<i} c_{ij} \be_j
\]
for some $\lambda_i \in \calO_0$ with $|\lambda_1| \geq \cdots \geq |\lambda_n|$
and $c_{ij} \in \Gsep$. Put $|\lambda_i| = p^{-\ell_i}$;
we may assume that $|\lambda_i| > |c_{ij}|$ for $j<i$.

Now suppose $\bv = \sum_{i=1}^n d_i \be_i$
is an eigenvector of $M$ of maximum slope. Equating the coefficients
of $\be_i$ in the equation $F\bv = \lambda_n \bv$ gives
\[
\lambda_n d_i = \lambda_i d_i^\sigma + \sum_{j>i} d_j^{\sigma} c_{ji}.
\]
Let $m$ be the smallest integer such that $|\lambda_m| = |\lambda_n|$;
by the hypothesis that $M$ is not isoclinic, we know that $m>1$.

We wish to show that $d_i \in \Gamma^{\alg(c)}$ for $i = 1,\dots,n$
by descending induction on $i$. For $i=n$, we have
$\lambda_n d_n = \lambda_n d_n^\sigma$, so $d_n \in \calO_0$.
For $m \leq i < n$, we have $d_i \in \Gsep$ by descending induction,
since $d_i = d_i^\sigma + x$ for $x \in \Gsep$.
Now suppose $i < m$;
assuming $d_{i+1}, \dots, d_n \in \Gamma^{\alg(c)}$, we have that
\begin{equation} \label{eq:filt1}
(\lambda_n/\lambda_i) d_i = d_i^\sigma + x \qquad (x \in \Gamma^{\alg(c)}, |x|<1).
\end{equation}
Let $|\lambda_n/\lambda_i| = p^{-e}$, so that $e \geq c$.
For $y \in \Galg$ and $\ell$ a nonnegative integer,
let $f_\ell(y)$ be the smallest real number $m$
such that there exists $z \in \Gsep$ with $|y^{\sigma^{\ell}}-z| < p^{-m}$.
Now $y \in \Gamma^{\alg(c)}$ if and only if $f_\ell(x) \geq c\ell$ for $\ell
\geq 0$.
Thus we have that $f_\ell(x) \geq c\ell$ for $m \geq 0$ and we wish to show
that $f_\ell(d_i) \geq c\ell$ for $\ell \geq 0$, by induction on $m$.
For $\ell=0$, we have $|x|< 1$ and so $|d_i|<1$, yielding $f_0(d_i)=0$.
As for $\ell>0$, from (\ref{eq:filt1}),
we have
\[
f_{\ell}(d_i) \geq
e+\min\{f_{\ell-1}(d_i), c\ell\}
\geq c\ell.
\]
We conclude that $d_i \in \Gamma^{\alg(c)}$.
\end{proof}

The part of the following lemma with $\phi$ mapping to $\Gamma$ is
\cite[Corollary~8.2]{bib:dej1}.
\begin{lemma} \label{slopefilt2}
Let $M$ be an $F$-crystal over $\Gcon$ admitting an injective map
$\map{\phi}{M}{\Gamma}$ (resp.\ $\map{\phi}{M}{\Gperfcon}$)
for which $\phi(F\bv) = \lambda \phi(\bv)^\sigma$
for all $\bv \in M$, with $\lambda \in \calO_0$.
Let $\ell$ be the highest (resp.\ lowest) slope of $M$. Then the
following results hold.
\begin{enumerate}
\item[(a)]
The multiplicity of $\ell$ as a slope of $M$ is $1$.
\item[(b)]
$|\lambda| = p^{-\ell}$.
\item[(c)]
$\phi^{-1}(\Gcon)$ is a rank $1$ submodule
of $M$ stable under $F$ with slope $\ell$.
\end{enumerate}
\end{lemma}
\begin{proof}
By Proposition~\ref{inj1}, the map $\map{\psi}{M \otimes_{\Gcon} \Gimmcon}{\Gimm}$
(resp.\ $\map{\psi}{M \otimes_{\Gcon} \Gsep}{\Galg}$)
induced by $\phi$ is also injective. Let $M_1$ be the first term of the descending slope
filtration of $M \otimes_{\Gcon} \Gimmcon$
(resp.\ the ascending slope filtration of $M \otimes_{\Gcon} \Gsep$).
We know $M_1$ is spanned by vectors $\bv$ with
$F\bv = \lambda_1 \bv$. On the other hand, $\psi(F\bv) = \lambda \psi(\bv)^\sigma$;
since we cannot have $\psi(\bv) = 0$ by injectivity, we must have $|\lambda| =
|\lambda_1| = p^{-\ell}$
and hence $\psi(\bv) \in \calO_0$. If $M_1$ has rank greater than 1, then some linear
combination of eigenvectors goes to 0, which is forbidden, so $M_1$ has rank 1
and we may choose a generator $\bv$ of $M_1$ with $\psi(\bv)=1$.

All that remains to be shown is that $\phi^{-1}(\Gcon) = M \cap \psi^{-1}(\Gcon)$
is nonempty.
Pick an isobasis $\be_1, \dots, \be_n$ of
$M$ over $\Gcon$, and write
$\bv = \sum_i c_i \be_i$. Now $\lambda = \sum_i c_i \phi(\be_i)$,
so by Proposition~\ref{inj1} again,
$\sum_i c_i \otimes \phi(\be_i) - \lambda \otimes 1 = 0$.
In particular, there must exist a nonzero linear combination of
the $\phi(\be_i)$ over $\Gcon$ which sums to 1. If $\sum d_i
\phi(\be_i) = 1$, then $\sum d_i \be_i \in \phi^{-1}(\Gcon)$,
and the latter is nonempty, as desired.
\end{proof}
It can be seen from the proof that the argument also applies to
an $F$-crystal $M$ over $\Gamma$ and a map $\map{\phi}{M}{\Gperf}$, 
the final conclusion then being that $\phi^{-1}(\Gamma)$ is
a rank 1 submodule of $M$ stable under $F$ with slope $\ell$.

\section{More on semistable crystals}

\begin{lemma}
Let $M$ an $F$-crystal over $\Gancon$ and $\be_1, \dots, \be_n$ an isobasis
of $M$ such that $F\be_i = \lambda_i \be_i$ with $\lambda_i \in \calO_0$. Then
any eigenvector of $M$ is an $\calO$-linear combination of
the $\be_i$.
\end{lemma}
\begin{proof}
Suppose $\bv \in M$ is such that $F\bv = \lambda \bv$, with $\lambda \in \calO$. Write
$\bv = \sum_i c_i$ with $c_i \in \Gancon$; then $\lambda c_i = \lambda_i
c_i^\sigma$. If $|\lambda| \neq |\lambda_i|$, then this equation implies
$c_i=0$; otherwise, we must have $c_i \in \calO$.
\end{proof}

Let $M$ be an $F$-crystal over $\Gcon$ under the standard Frobenius $\sigma_t$.
Recall that $M$ is defined to be \emph{semistable}
if and only if it isomorphic to $M' \otimes_\Omega \Gcon$ for some
$F$-crystal $M'$ over $\Omega$.

\begin{lemma} \label{nocon}
If $M$ is an $F$-crystal over $\Omega$ and $N$ is an
$F$-stable submodule (but not necessarily a subcrystal)
of $M \otimes_\Omega \Gcon$, then there exists a
sub $F$-crystal $N'$ of $M$
such that $N = N' \otimes_\Omega \Gcon$.
\end{lemma}
\begin{proof}
By replacing $M$ with a suitable exterior power,
we may reduce to the case where $N$ has rank 1. Let $\bv$ be a
generator of $N$, so that $F\bv = c \bv$ for some $c \in \Gcon$.  By
Dwork's trick \cite[Lemma~4.3]{bib:me3},
we may (after enlarging $\calO$) choose a basis $\be_1,
\dots, \be_n$ for $M$ over $\Oan$ such that $F\be_i = \lambda_i \be_i$
with $\lambda_i \in \calO_0$. Write $\bv = \sum c_i \be_i$ with $c_i \in
\Gancon$; then we have $\lambda_i c_i^\sigma = c c_i$ for $i=1,
\dots, n$.  In particular, $(c_i/c_j)^\sigma = (\lambda_i/
\lambda_j) c_i/c_j$ whenever $c_j \neq 0$.  From this we deduce that
for any $i,j$ such that $c_i, c_j \neq 0$, $|\lambda_i|=
|\lambda_j|$ (since the equation $x^\sigma = \lambda x$ has no solutions in $\Gancon$ for
$|\lambda| \neq 1$). Moreover, any two nonzero $c_i$ are multiples of one another by elements
of $\calO_0$ (since these are the only solutions of $x^\sigma =x$ in $\Gancon$).
Therefore
we can write $\bv = d \bw$, where $\bw$ is an eigenvector of $M$ over $\Oan$
and $d \in \Gancon$.

By \cite[Corollary~4.10]{bib:me3}, we may write $d=ef$
with $e \in \Gcon[\fp]$ and $f \in \Oan^*$;
we may
shift factors of $\fp$ to
or from $f$ to ensure that $e \in \Gcon^*$. Now $e \bv =
f^{-1} \bw$ is defined over $\Oan \cap \Gcon = \Omega$ and is a
direct summand of $M$ over $\Gancon$, hence also over $M$.
Additionally, $F(e \bv)= k e \bv$ with $k = f^\sigma \lambda/f \in
\Oan^* \cap \Gcon = \Omega$.

Finally, note that since the span of $e \bv$ is a direct summand over $\Omega$
up to isogeny,
$F$ acts on it through a matrix invertible
over $\Omega[\fp]$. Thus this span is actually a subcrystal, as
desired.
\end{proof}

\begin{cor} \label{cor:exact}
Let $M$ be an $F$-crystal over $\Gcon$ and
\[
0 \to M_1 \to M \to M_2 \to 0
\]
an exact sequence of $F$-crystals. If $M$ is semistable,
then $M_1$ and $M_2$ are
semistable.
\end{cor}
\begin{proof}
The fact that $M_1$ is semistable follows from Lemma~\ref{nocon}.
Now regarding $M_1$ and $M$ as actually being defined over $\Omega$,
we can construct a basis for $M_1$ which extends to a basis for $M$
(since the action of $F$ on $M_1$ is invertible up to scalars). The
quotient is defined over $\Omega$ and is isomorphic to $M_2$.
\end{proof}

The converse of this corollary need not hold in general: suppose $M$ is
a rank 2 crystal whose $F$-action on some basis is
$\begin{pmatrix} 1 & t^{-1} \\ 0 & p \end{pmatrix}$. Then $M$ is not semistable
even though it has a semistable rank 1 submodule with semistable quotient.

The above example also shows that over $\Gcon$ or even
$\Gancon$, exact sequences of $F$-crystals
may fail to be split. On the other hand, exact sequences do split under
certain circumstances, such as the following.
\begin{prop} \label{prop:johan71}
Let $0 \to M_1 \to M \to M_2 \to 0$ be an exact sequence of
semistable $F$-crystals
over $\Gcon$, and assume the slopes of $M_1$ exceed those of $M_2$.
Then the exact
sequence splits over $\Gcon$.
\end{prop}
\begin{proof}
By taking a suitable exterior power, performing an isogeny and
dividing the action
of $F$ by a scalar, we reduce to the case where $M_1$ has dimension 1.
Then the statement follows immediately from the previous lemma and from
\cite[Proposition~7.1]{bib:dej1}.
\end{proof}

\section{Descent of morphisms}

For convenience, we recall the statement of Tsuzuki's conjecture
\cite[Conjecture~2.3.2]{bib:tsu4},
recast in the notations of this paper.
\begin{conj} \label{main2}
Any morphism over $\Gamma$ between $(F, \nabla)$-crystals over $\Gcon$
is obtained from a morphism over $\Gcon$ by extension of scalars.
\end{conj}
Equivalently, for $M$ an
$(F, \nabla)$-crystal over $\Gcon$, any eigenvector in $M \otimes_{\Gcon}
\Gamma$ is actually in $M$. (A morphism 
from $M$ to $M'$ can be viewed as an eigenvector in
$\Hom(M, M') = M^* \otimes M'$, where $M^*$ is the dual crystal.)

We first prove the extension of \cite[Theorem~9.1]{bib:dej1}
to the logarithmic case;
the proof is obtained from de~Jong's by reading closely and eliminating
all references to the connection.

\begin{theorem} \label{johan}
Any morphism over $\Gamma$ between $F$-crystals over $\Omega$
is obtained from a morphism over $\Omega$ by extension of scalars.
\end{theorem}
\begin{proof}
Equivalently, we show that for $M$ an $F$-crystal over $\Omega$, and
$\bv \in M \otimes_\Omega \Gamma$ with $F\bv = \lambda \bv$
for some $\lambda \in \calO_0$, that $\bv \in M$.
Form the dual $M^* = M^*(\ell)$ for some large integer $\ell$;
from $\bv$ we obtain a map $\map{\phi}{M^*}{\Gamma}$
such that $\phi(F\bx) = \lambda \phi(\bx)^\sigma$ for
all $\bx \in M^*$. By Lemma~\ref{nocon}, there exists an $F$-stable
submodule $N$ of $M^*$ whose extension to $\Gcon$ is the kernel
of the map $M^* \otimes \Gcon$ to $\Gamma$. Let $M_1 = M^*/N$.

By Lemma~\ref{slopefilt2}, the highest slope of
$M_1 \otimes \Gcon$ is $\ell + v_p(\lambda)$,
has multiplicity 1, and has eigenspace $M_2 = \phi^{-1}(\Gcon)$.
By Corollary~\ref{cor:exact}, $M_3 = M_1/M_2$ is semistable,
and the slopes of $M_3$ are all less than $\ell+v_p(\lambda)$,
so Proposition~\ref{prop:johan71}
gives a direct sum decomposition $M_1 = M_2 \oplus M_3$. However,
Lemma~\ref{slopefilt} forces $M_3$ to map to 0 under $\phi$. Since
$\phi$ is injective on $M_1$, $M_3 = 0$. Thus $M_1$ is one-dimensional and
constant over $\Omega$, and $\phi(M^*) = \phi(M_1) = \Omega$, as desired.
\end{proof}

As mentioned earlier, we are currently unable to resolve Tsuzuki's conjecture
in full generality. However, using Tsuzuki's theorem that isoclinic crystals
are potentially semistable, we can resolve the conjecture for morphisms
between crystals with at most two distinct slopes each. In particular,
this includes the crystals arising from the cohomology of
ordinary abelian varieties.

\begin{cor} \label{cor:twoslope}
Conjecture~\ref{main2} holds for a morphism between crystals with at most
two distinct slopes each. In particular, it holds for crystals of rank $2$.
\end{cor}
\begin{proof}
We first show that if $M$ is an $(F, \nabla)$-crystal all of whose slopes
are equal to one of $a,b,c$, where $a<b<c$, and $\bv$ is
an eigenvector $\bv$ over $\Gamma$ with slope $b$ and $\nabla \bv = 0$,
then $\bv$ is defined over
$\Gcon$. Construct the dual $M^* = M^*(\ell)$ of $M$, with
slopes $\ell-c, \ell-b, \ell-a$.
Let $\phi: M^* \to \Gamma$ be the map induced by $\bv$,
let $M_1$ be the kernel of $M^*$, let $M_2 = M^*/M_1$,
%WATCH OUT!
and let $M_3$ be the top eigenspace of $M_2$,
which has rank 1. Then $M_2/M_3$ is isoclinic of slope $\ell-c$, and so is
potentially semistable, in fact potentially constant,
by \cite[Theorem~5.1.1]{bib:tsu1}.

Take a finite extension $\GL_{\con}$ of $\Gcon$
over which $M_2/M_3$ is constant.
Take $\bv \in M_3$ with $F\bv = \lambda \bv$ for some $\lambda \in \calO_0$;
then we can extend
$\bv$ to a basis $\bv, \bx_1, \dots, \bx_n$ of $M_2$ with $F\bx_i
= \mu \bx_i + c_i \bv$ for some $\mu \in \calO_0$ with $|\lambda| < |\mu|$.
Applying $\nabla$ to this equation gives
\[
pF(\nabla \bx_i) - \mu \nabla \bx_i = c_i' \bv.
\]
This equation uniquely determines $\nabla \bx_i$: it must have the form
$b\bv$, where $y \in \GLcon$ satisfies
and $p\lambda b^\sigma - \mu b = c_i'$. On the other hand, the equation
$\lambda d^\sigma - \mu d = c_i$ has a unique solution with $d \in \GL$,
and differentiating gives $p\lambda (d')^\sigma - \mu d' = c_i'$. Thus
we must have $d' = b$ and so $d \in \GLcon$, whence $\bx_i - d \bv$ is
an eigenvector of $M_3$. The upshot is that $M_2$ has a complement in $M_3$
over $\GLcon$, so
$\phi$ maps into $\GLcon$ and $\bv$ is defined over $\GLcon \cap
\Gamma = \Gcon$.

Suppose we are given a morphism $f: M_1 \to M_2$
over $\Gamma$ between two $(F, \nabla)$-crystals $M_1$ and $M_2$
over $\Gcon$ with at most two distinct slopes each.
If the slopes of $M_1$ and $M_2$ are the same, then $M_1^* \otimes M_2$ has
three distinct slopes and $f$ corresponds to an eigenvector of the middle
slope; by the result of the previous paragraph, $f$ is
then defined over $\Gcon$.

Now suppose the slopes of $M_2$ are not the same as the slopes of $M_1$.
Then the image
of $f$ must be an isoclinic sub-$(F,\nabla)$-crystal of $M_2$;
again by \cite[Theorem~5.1.1]{bib:tsu1}, the eigenvectors in the subcrystal
are defined over
a finite separable extension $\GL_{\con}$ of $\Gcon$. Thus over $\GL_{\con}$,
$f$ maps $M_1$ to a direct sum of trivial crystals.
Composing with the projection
onto one factor gives a map from $M_1$ to one trivial crystal,
which corresponds
to an eigenvector of $M_1^*$ over $\GL$. By the first claim above,
this eigenvector is defined over $\GL_{\con}$;
after repeating for each
factor in the direct sum, we conclude that $f$ is defined over $\Gamma
\cap \GL_{\con} = \Gcon$.
\end{proof}

\section{Trivial representations}

Theorem~\ref{johan} can be used to answer
a question of Katz \cite[p.\ 162-163]{bib:katz}. Namely,
if $M$ is a unipotent $F$-crystal over $\Gcon$, then the representation of
$\Gal(\Galg/\Gperf)$ on
the eigenvectors over $\Galg$ is trivial; for $M$ over $\Omega$,
Katz asked whether conversely the triviality of the Galois
representation implies that $M$ is unipotent over $\Omega$.
To see that this is so, note that the triviality
of the representation implies that the eigenvectors of $M$ are all defined
over $\Gperf$. In particular, the eigenvectors of lowest slope are defined
over $\Gperf \cap \Gsep = \Gamma$ by Lemma~\ref{slopefilt1a}. Thus by
Theorem~\ref{johan}, these are defined over $\Omega$,
so we may quotient by these and repeat the argument.

One may also ask whether an $(F, \nabla)$-crystal over
$\Gcon$ is unipotent as an $F$-crystal
if and only if its corresponding Galois representation
is trivial. We now show that this is the case.

\begin{theorem}
Let $M$ be an $(F, \nabla)$-crystal over $\Gcon$. If the Galois
representation associated to $M$ is trivial, i.e., if $M$ becomes
constant over $\Gperf$, then $M$ is unipotent as an $F$-crystal.
\end{theorem}
\begin{proof}
We induct on the dimension of $M$.
Suppose $M$ becomes constant over $\Gperf$. Let $c$ be the highest slope
of $M$, and let $M^* = M^*(\ell)$ be a dual of $M$.
Then the eigenvectors of $M^*$
of lowest slope $\ell-c$
are defined over $\Gperf \cap \Gsep = \Gamma$. Let $\bv$
be one of these eigenvectors, with $F\bv = \mu \bv$;
then $\bv$ corresponds to a map $f:
M \to \Gamma$ such that $f(F\bw) = p^\ell \mu^{-1} f(\bw)^\sigma$.
By Lemma~\ref{slopefilt2}, 
$M/\ker(f)$ has highest slope $c$ with multiplicity 1 and eigenvector
defined over $\Gcon$.
By the induction hypothesis, $\ker(f)$ is unipotent over $\Gcon$,
as then is $M$.
\end{proof}

\section{Geometric crystals}
\label{sec:geometric}

In this section, we establish an assertion to the effect that
``crystals coming from geometry are overconvergent and
potentially semistable''. Namely,
let $\psi: V \to \Spec k((t))$ be a smooth proper morphism. (By the definition
of properness, $\psi$ is automatically of finite type.)
Then Ogus \cite[Section~3]{bib:ogus} showed how to define the derived
functor $R^q f_* \calO_{V, \cris}$ (modulo torsion)
as a convergent $F$-isocrystal over $k((t))$, which is to say,
as an $(F, \nabla)$-crystal over $\Gamma$. (Note: here and throughout
this section, $\calO$ indicates a structure sheaf, not a finite
extension of $W(k)$.) We wish to
show that this crystal is in fact overconvergent.

\begin{theorem} \label{geom}
Let $\psi: V \to \Spec k((t))$ be a smooth proper morphism,
and for $i \geq 0$, let $M_i = R^i f_* \calO_{V,\cris}$ as an $(F, \nabla)$
crystal over $\Gamma$.
Then $M_i$ is overconvergent (i.e.\ descends to a crystal over $\Gcon$)
and potentially semistable.
\end{theorem}
\begin{cor}
Morphisms of $(F, \nabla)$-crystals arising as $R^i f_* \calO_{V, \cris}$
(or more generally, arising from motives) descend from $\Gamma$ to
$\Gcon$.
\end{cor}

The argument is motivated by Berthelot's proof \cite{bib:ber1}
of the finite dimensionality
of the rigid cohomology of an arbitrary variety (with constant coefficients).
Presumably this result can also be demonstrated for $V$ not
necessarily smooth, using rigid cohomology in place of 
crystalline cohomology. The results of \cite{bib:ber2}
will no doubt facilitate this task, but we have not yet seen this paper.
Additional possible refinements would be to show that the structure sheaf on $V$
can be replaced by an arbitrary crystal of modules of finite type,
perhaps even with $V$ not proper as long as the crystals themselves are required
to be overconvergent in a suitable sense.

We first describe why Theorem~\ref{geom} holds in case $V$ admits a ``semistable''
extension across $\Spec k[[t]]$. Namely, suppose there exists a proper
morphism $\overline{\psi}: \overline{V} \to \Spec k[[t]]$ which
extends $\psi$, such that $\overline{V}$ is a regular scheme and
the special fiber of $\overline{\psi}$ is reduced with strict normal
crossings. Then one gets a smooth map $\overline{\psi}:
(\overline{V}, A) \to (\Spec k[[t]], B)$ of schemes with fine logarithmic
structures, in the sense of Kato \cite{bib:kato1}.

The construction of Hyodo-Kato \cite[Proposition~2.24]{bib:hk} shows that
$R^i \overline{\psi}_* \calO^{\log}_{\overline{V}, \cris}$ (again
modulo torsion) is an $(F, \nabla)$-crystal over $\Omega$. More specifically,
one uses $\Omega$ as a test object in the crystalline site over $k[[t]]$
and obtains a locally free module over $\Omega$ with actions of $F$ and 
$\nabla$ from the Frobenius automorphism and Gauss-Manin connection,
respectively,
of $V$. The result of Hyodo-Kato shows that the kernel and cokernel of $F$
are finite, so that the result is indeed an $(F, \nabla)$-crystal.
(This process
is a logarithmic analogue of the argument of \cite{bib:ber3} reducing
the extension theorem for $p$-divisible groups in equal characteristic $p$ to
a statement about $F$-crystals over $\Omega$.)

To prove Theorem~\ref{geom}, we must work around the fact that it is not
known whether $\psi$ must extend to a log-smooth morphism over $k[[t]]$,
even after a finite base extension. As is now customary, the workaround
involves de~Jong's results on alterations, specifically his
semistable alterations theorem (a relative form of weak
resolution of singularities). The result in \cite{bib:dej0} is actually somewhat more
general; we cite the precise formulation we will be using.
\begin{theorem}[\cite{bib:dej0}, Theorem~6.5] \label{thm:alter}
Let $X$ be a scheme of finite type over a trait (spectrum of a complete
discrete valuation ring) $S$, whose generic fiber is reduced and geometrically
irreducible, 
and $Z \subset X$ a proper
closed subset containing the special fiber of $X$. Then there
exists a trait $S_1$ finite over $S$, a proper variety $X_1$ over $S_1$,
an alteration $\phi: X_1 \to X$ over $S$ and an open immersion $j_1:
X_1 \to \overline{X_1}$ of varieties over $S_1$, with the following
properties.
\begin{enumerate}
\item $\overline{X_1}$ is a projective variety over $S_1$ whose generic fiber
is geometrically irreducible.
\item The pair $(\overline{X_1}, \phi^{-1}(Z))$ is strictly semistable. 
In particular, $\phi^{-1}(Z)$ is a strict normal crossings divisor
on $\overline{X_1}$. (For the full definition, see \cite[Section~6.3]{bib:dej0}.)
\end{enumerate}
\end{theorem}
%DIAGRAM.

We will also need a quick lemma on the action of an alteration on cohomology.
\begin{lemma}
Let $\map{\phi}{X}{Y}$ be a surjective map between smooth, irreducible,
proper varieties $X$ and $Y$ of the same dimension. Then
$\map{\phi^*}{H^i_{\cris}(Y)}{H^i_{\cris}(X)}$ is injective. Moreover,
there exists a projector mapping $H^i_{\cris}(X)$ to the image of $\phi^*$.
\end{lemma}
\begin{proof}
Define the pushforward map $\phi_*$ by applying Poincar\'e duality
$H^{i}_{\cris}(V)^* \cong H^{\dim V - i}_{\cris}(V)$ to the transpose
of $\phi^*$.

Let $1_Y$ denote the fundamental class on $Y$.
For arbitrary classes $\omega$ and $\eta$ on $Y$ of complementary
dimension, we have
\begin{align*}
\phi_* \phi^* \omega \cup \eta &= \phi_* (\phi^* \omega \cup \phi^* \eta) \\
&= \phi_* \phi^* (\omega \cup \eta) \\
&= \phi_* \phi^* 1_Y \cup \omega \cup \eta.
\end{align*}
On the other hand, $\phi^* 1_Y = 1_X$, the fundamental class on $X$. On the other hand,
$1_X$ is Poincar\'e dual to the cycle class of a point, whose pullback is Poincar\'e
dual to the class of the preimage of that point (assuming the point is chosen in the
dense open set on which $\phi$ is finite). Thus $\phi_* 1_X = d 1_Y$, where
$d$ is the generic degree of $\phi$. In particular, $d \neq 0$, so $\phi_* \phi^*$
is injective, as then must be $\phi^*$, and the desired projector is
$\frac{1}{d} \phi^* \phi_*$.
(For the relevant inputs into this proof related to crystalline cohomology,
see \cite{bib:ber0}, specifically Section~VI.3.3 for
cycle classes and Section~VII.2.2 for Poincar\'e duality.)
\end{proof}

\begin{proof}[Proof of Theorem~\ref{geom}]
We first consider the case where $V$ is projective. In this case, we may
extend $V$ to a scheme $X$ of finite type
over $\Spec k[[t]]$ by choosing
an embedding of $V$ in a projective space over $k((t))$ and taking its
Zariski closure in $\PP^n_{k[[t]]}$.

Now let $S_1, X_1, \phi, j_1$ be as in Theorem~\ref{thm:alter}.
Then $R_i\psi_* \calO_{X_1,\cris}$ is isomorphic
over $\Gamma$ to the log-crystal $M_{i,1} =
R_i g_* \calO_{\overline{X_1},\cris}^{\log}$,
hence is semistable.

Now the surjection $(X_1)_{\eta_1} \to V \times_{\eta} \eta_1$
induces a map
$M_i \otimes_{\Gamma} \Gamma^{k((u))} \to R_i \psi_* \calO_{W,\cris}$ of
crystals over $k((u))$, which gives an injection $M_i \otimes_\Gamma 
\Gamma_{k((u))} \hookrightarrow M_{i,1}$
of $(F, \nabla)$-crystals over $\Gamma$ by the previous lemma. Moreover, the projector
constructed in the proof of the lemma is an endomorphism of $M_{i,1}$
over $\Gamma$; by Theorem~\ref{johan}, this endomorphism is actually defined
over $\Omega$, as then is its image
$M_i \otimes \Gamma^{k((u))}$. Thus $M_i \otimes \Gamma^{k((u))}$ is semistable,
and in particular is also overconvergent; since $\Gamma^{k((u))}_{\con}$ is
a finite extension of $\Gamma^{k((t))}_{\con}$,
we conclude that $M_i$ is also overconvergent, and is potentially semistable.

To handle the case of $V$ arbitrary, recall that by Chow's Lemma, there exists
$V_1$ projective such that $V_1 \to V$ is a surjective birational morphism. By the same
argument as in the previous paragraph, but applied to the map
$R_i \psi_* \calO_{V,\cris} \to R_i \psi_* \calO_{V_1,\cris}$,
we conclude that the latter being potentially semistable implies the same for the 
former.
\end{proof}

\section*{Acknowledgments}
This work is based on the author's doctoral dissertation, written at
the Massachusetts Institute of Technology under the supervision of
Johan de~Jong. The author was supported by
a National Science Foundation Postdoctoral
Fellowship.

%% Some of these refs are extraneous.

\end{document}